\documentclass{amsart}

\usepackage{amsthm, amssymb, amsmath}
\usepackage{mathtools}
\usepackage{xcolor}
\mathtoolsset{showonlyrefs}

\newtheorem{theorem}{Theorem}
\newtheorem{corollary}[theorem]{Corollary}

\newtheorem{definition}[theorem]{Definition}
\newtheorem{lemma}[theorem]{Lemma}

\def \R {\mathbb R}
\def \C {\mathbb C}
\def \N {\mathbb N}
\newcommand{\vol}[1] {|#1|_n}
\def \B{\mathbb B}

\def \S {\mathcal S}
\def \I {\mathcal I_n}
\newcommand {\D}[1][q] {(-\Delta)^{{#1}/2}}
\def \Rad {{\mathcal R}}
\newcommand {\Radq}[1] {\frac 1{\cos(\pi q/2)}\mathcal R #1(\xi, \cdot)_t^{(q)}(0)}
\def \dovr {d_{\operatorname{ovr}}}

\title[]{Radon transforms with small derivatives and distance inequalities for convex bodies}
\author{Juli\'an Haddad}
\address{Departamento de An\'alisis Matem\'atico, Universidad de Sevilla, Sevilla, Espa\~na}
\email{jhaddad@us.es}

\author{Alexander Koldobsky}
\address{Department of Mathematics, University of Missouri-Columbia, Columbia MO 65211, USA}
\email{koldobskiya@missouri.edu}

\begin{document}
\maketitle
\begin{abstract} Generalizing the slicing inequality for functions on convex bodies from \cite{bib_koldobsky_slicing},
it was proved in \cite{bib_gregory_fractional_lower_bound} that there exists an absolute constant $c$ so that for any $n\in \N$,
any $q\in [0,n-1)$ which is not an odd integer, any origin-symmetric convex body $K$ of volume one in $\R^n$ and any infinitely 
smooth probability density $f$ on $K$ we have 
$$\max_{\xi \in S^{n-1}} \Radq f \ge \left( \frac {c(q+1)}{n}\right)^{\frac{q+1}2}.$$
Here $\mathcal Rf(\xi,t)$ is the Radon transform of $f$, and the fractional derivative of the order $q$ is taken with respect to the variable $t\in \R$
with fixed $\xi\in S^{n-1}.$ In this note we show that there exist an origin-symmetric convex body $K$ of volume 1 in $\R^n$ and a continuous probability density
$g$ on $K$ so that
$$\max_{\xi\in S^{n-1}} \Radq g \leq \frac 1{\sqrt n} (c(q+1))^{\frac{q+1}2}.$$
In the case $q=0$ this was proved in \cite{bib_klartag_lower_bound_with_loglog, bib_livschyts_lower_bound}, and it was  used there to obtain a lower estimate for the maximal outer volume ratio distance from an arbitrary origin-symmetric convex body $K$ to the class of intersection bodies.
We extend the latter result to the class $L_{-1-q}^n$ of bodies in $\R^n$  that embed in $L_{-1-q}.$ Namely, for every $q\in [0,n)$ there exists an origin-symmetric convex body $K$ in $\R^n$ so that 
 $\dovr(K, L_{-1-q}^n) \ge cn^{\frac 1{2(q+1)}}.$
\end{abstract}

\section{Introduction}

The slicing problem for functions asks for the greatest constant $s_n,$ depending only on the dimension $n,$ such that for every origin-symmetric convex body $K$ in $\R^n$ of volume 1 and any continuous probability density $f$ supported in $K$ there exists a hyperplane $H$ passing through the origin for which
\[\int_H f \geq s_n.\]

This problem was put forward in \cite{bib_koldobsky_hyperplane_inequality, bib_kold_sqrt_estimate, bib_koldobsky_slicing} and further investigated in \cite{CGL, bib_bobkov_moments}. The case where $f \equiv 1$ is the slicing problem of Bourgain, for which the best lower bound so far, due to Klartag  \cite{bib_klartag_slicing}, is $c (\log n)^{-1/2}$.

More precisely, let us define
\[s_n = \inf_{f,K} \max_{\xi \in S^{n-1}} \int_{K\cap \xi^\perp} f,\]
where the infimum is taken over the class of all even continuous probability densities defined on symmetric convex bodies $K$ of volume $1$ in $\R^n$.
Here $\xi^\perp=\{x\in \R^n: (x,\xi)=0\}$ is the central hyperplane perpendicular to $\xi.$

It was proved in \cite{bib_kold_sqrt_estimate} that $s_n \geq \frac 12 n^{-1/2}.$ A more general result was proved in \cite{bib_koldobsky_slicing},
as follows.

\begin{theorem}
	\label{res_upper_bound_sqrt}
	Let $K$ be an origin-symmetric convex body in $\R^n$, and let $f$ be an even continuous non-negative function on $K$. Then
	\[
		\int_K f \leq \dovr(K,\I) \frac n{n-1} c_n  \vol{K}^{1/n} \max_{\xi\in S^{n-1}}\int_{K \cap \xi^\perp} f, 
	\]
	where $c_n = \frac{\vol{B_2^n}^{\frac{n-1}n}}{|B_2^{n-1}|_{n-1}} < 1$
	and $\dovr(K,\I)$ is the outer volume ratio distance from $K$ to the class of intersection bodies $\mathcal{I}_n$ (see Section \ref{sec_preliminaries} for the definition).
\end{theorem}
The outer volume ratio distance to a class of convex bodies $\mathcal C$ is defined by
\[\dovr(K, \mathcal C) = \inf \left\{ \left( \frac{\vol{D}}{\vol{K}}\right)^{1/n}: K\subseteq D, D \in \mathcal C \right\}.\]
Since the class of intersection bodies contains ellipsoids and $K$ is symmetric, by John's theorem one always has $\dovr(K,\I) \leq \sqrt n$.
The constant $\dovr(K,\I)$ may be smaller for special classes of bodies $K$. For example,
if $K$ is an unconditional convex body then $\dovr(K,\I) \leq e$ (see \cite{bib_koldobsky_slicing}), and if $K$ is the unit ball of a subspace of $L_p([0,1])$ then $\dovr(K,\I) \leq \sqrt p$ (see \cite{M, bib_koldobsky_lp}).

On the other hand, the upper bound $s_n \leq c n^{-1/2}$ was established in \cite{bib_klartag_lower_bound_with_loglog} up to a logarithmic term, and the proof was refined in 
 \cite{bib_livschyts_lower_bound} to remove the logarithmic factor.
\begin{theorem}
	\label{res_example_KL}
There exists a symmetric convex body $K \subseteq \R^n$ and an even continuous function $f:\R^n \to \R$ for which
	\[\int_K f \geq c, \qquad \int_{K \cap \xi^\perp} f \leq \frac c {\sqrt n} \text{ and } \qquad \vol{K} \leq c^n\]
	for every $\xi \in S^{n-1}$, where $c>0$ is an absolute constant.
	Moreover, $f$ belongs to the Schwartz space, and we have $\sqrt n \B \subseteq K$.
\end{theorem}

It immediately follows from Theorems \ref{res_upper_bound_sqrt} and \ref{res_example_KL} that 
$$c\sqrt{n} \le \sup_K \dovr(K,\I) \le \sqrt n,$$ 
where $c>0$ is a universal constant, and the supremum is taken over all origin-symmetric convex bodies $K$ in $\R^n.$

By applying a similar strategy to moments of functions on convex sets, it was shown in \cite{bib_bobkov_moments} (lower bound) and \cite{KPZ} (upper bound) that
$$c_1\sqrt{\frac np}\le \sup_K \dovr(K, L_p^n) \le c_2 \sqrt{\frac np},$$
where $c_1,c_2$ are absolute constants, and $L_p^n$ is the class of unit balls of $n$-dimensional spaces that embed in $L_p([0,1]),\ p\ge 1$.
\smallbreak

For a continuous function $f:\R^n \to \R$, the Radon transform of $f$ is defined by 
\[\Rad f(\xi, t) = \int_{\xi^\perp+t \xi} f, \qquad \xi\in S^{n-1},\ t\in \R.\]
Clearly, Theorems \ref{res_upper_bound_sqrt} and \ref{res_example_KL} can be stated as bounds for $\Rad f(\xi, 0)$, and
\[s_n = \inf_{f,K} \max_{\xi \in S^{n-1}} \Rad (\chi_K f)(\xi, 0)\]
where $\chi_K$ is the indicator function of $K$.

We denote by $\mathcal R f (\xi, \cdot)_t^{(q)}(0)$ the fractional derivative of order $q$ of the function $t \mapsto \Rad f(\xi, t)$ at $t=0$ (see Section \ref{sec_preliminaries} for the definition).

The following Theorem is a variant of \cite[Theorem 2]{bib_gregory_fractional_lower_bound}.

\begin{theorem}
	\label{res_upper_bound_fractional}
	Suppose $K \subseteq \R^n$ is an infinitely smooth origin-symmetric convex body, $f:\R^n \to \R$ is a non-negative even infinitely smooth function and $q\in(-1,n-1)$ is not an odd integer. Then
	\begin{multline}
		\int_K f \leq \frac n{(n-q-1)2^q \pi^{\frac{q-1}2}\Gamma(\frac{q+1}2) }\vol{K}^{\frac{q+1}n} \\ \times \dovr(K,L_{-1-q}^n)^{q+1} \max_{\xi \in S^{n-1}} \Radq f
	\end{multline}
	where $L_{-1-q}^n$ is the class of star bodies in $\R^n$  that embed in $L_{-1-q}$ (see Section \ref{sec_preliminaries}).
\end{theorem}
Note that Theorem \ref{res_upper_bound_fractional}, as stated in \cite[Theorem 2]{bib_gregory_fractional_lower_bound} requires $f$ to have support inside $K$.
But a quick inspection of the proof shows that this is actually not necessary, and the proof carries over without modifications.



In this note we adapt the example of Theorem \ref{res_example_KL} to obtain a function with small fractional derivatives of the Radon transform, thus complementing Theorem \ref{res_upper_bound_fractional}.
\begin{theorem}
	\label{res_fractional_example} Let $q\in [0,n-1),$ $q$ is not an odd integer.
	There exists an infinitely smooth, symmetric convex body $D \subseteq \R^n$ of volume $1$, and an even continuous function $g:\R^n \to \R$ for which $\int_D g = 1$ and
	\[\max_{\xi\in S^{n-1}} \Radq g \leq \frac 1{\sqrt n} (c(q+1))^{\frac{q+1}2}, \qquad \forall \xi\in S^{n-1}.\]
\end{theorem}

Note that the body $D$ does not depend on $q$.
\smallbreak
The reader should compare Theorem \ref{res_fractional_example} with \cite[Formula (7)]{bib_gregory_fractional_lower_bound}. Namely, for every $f,K$ under the conditions of Theorem \ref{res_upper_bound_fractional},
\begin{equation}
	\label{eq_upper_bound_fractional_normalized}
	\max_{\xi\in S^{n-1}} \Radq f \geq \frac{(c(q+1))^{\frac {q+1}2}}{(\dovr(K,L_{-1-q}^n))^{q+1}} \ge  \left( \frac {c(q+1)}{n}\right)^{\frac{q+1}2}.
\end{equation}

When $q=2k-1,\ k\in \N,$ we have
\begin{multline}
\lim_{q\to 2k-1} \Radq f \\
= (-1)^{k}(2k-1)!  \int_0^\infty t^{-2k}\left(Rf(\xi,t)-\sum_{j=0}^{k-1} \frac{t^{2j}}{(2j)!}(Rf(\xi,t))_t^{(2j)}(0)\right) dt,
\end{multline}
and Theorems \ref{res_upper_bound_fractional} and \ref{res_fractional_example} hold with this quantity.
\smallbreak

As a consequence of Theorem \ref{res_fractional_example} we obtain
\begin{corollary}
	\label{res_distance_bound}
	For $q \in [0,n-1)$ there exists an origin-symmetric convex body $K$ in $\R^n$ such that
	\[ \dovr(K,L_{-1-q}^n) \geq c_1 n^{\frac 1{2(q+1)}}.  \]
	\end{corollary}
	Combining this estimate with the upper estimate from \cite[Propositions 2 and 3]{bib_gregory_fractional_lower_bound}, we get
	$$c_1 n^{\frac 1{2(q+1)}}\le \sup_K \dovr(K,L_{-1-q}^n)\le \min\left\{c_2\sqrt{\frac{n\log^3(\frac{ne}{q+1})}{q+1}},\sqrt{n}\right\},$$
	where $c_1,c_2$ are absolute constants, and the supremum is taken over all origin-symmetric convex bodies in $\R^n.$	

\section{Notation and Auxiliary Results}
\label{sec_preliminaries}

A star body is a set of the form
\[K =\{0\} \cup \{x \in \R^n\setminus \{0\}: |x|_2 \leq \rho_K(x/|x|_2)\}\]
where $\rho_K:\S^{n-1} \to (0,\infty)$ is continuous and bounded.
The radial function $\rho_K$ is extended to $\R^n \setminus \{0\}$ as an homogeneous function of degree $-1$, this is $\rho_K(\lambda x) = \lambda^{-1} \rho_K(x)$.
We denote the volume of $K$ in $\R^n$ by $\vol{K}.$

The Minkowsky functional of a star body $K$ is $\|x\|_K = \rho_K(x)^{-1}$.
A convex star body $K$ is called a convex body, and $K$ is said to be infinitely smooth if $\rho_K$ is an infinitely smooth function of the sphere.

The radial metric in the class of star bodies is defined by
\[d_\rho(K,L) = \sup_{\xi \in S^{n-1}} |\rho_K(\xi) - \rho_L(\xi)|.\]

We use the Fourier transform of distributions and refer the reader to \cite{bib_kold_book} for details.
For $f \in \S$, the Schwartz space of rapidly decreasing infinitely differentiable functions, the Fourier transform is defined by
\[\hat f(\xi) = \int_{\R^n} f(x) e^{-i (x,\xi)} dx.\]
The Fourier transform of a distribution $f$ is defined by $\langle \hat f, \varphi \rangle = \langle f, \hat \varphi \rangle$ for every $\varphi \in \S$.
With this normalization we have $\hat{\hat f} = (2\pi)^n f$ for every even distribution $f$.

A star body $K$ in $\R^n$ embeds (isometrically) in $L_p,\ p>0$ if there exists a finite Borel measure $\mu$ on $S^{n-1}$ so that
$$\|x\|_K^p= \int_{S^{n-1}} |(x,\xi)|^p d\mu(\xi),\qquad \forall x\in \R^n.$$
We denote by $L_p^n$ the set of all such bodies.

Following Lutwak \cite{L}, we define the intersection body of a star body $K$ as the star body $IK$ for which $\rho_{IK}(\xi) = |K \cap \xi^\perp|_{n-1}$.
A more general class of bodies was introduced in \cite{Ko}.
Given two star bodies $K,L \subseteq \R^n$ and an integer $k \in [1, n]$ we say that $K$ is the $k$-intersection body of $L$ if
\[|K \cap H^\perp|_{k} = |L \cap H|_{n-k},
\]
for all $n-k$ dimensional subspaces of $\R^n$.
The class of $k$-intersection bodies is defined as the closure in the radial metric of $k$-intersection bodies of star bodies.
For $k=1$ we get the class of intersection bodies ${\mathcal I}_n.$

We say that a star body $K$ embeds in $L_{-p}$ if $(\|\cdot\|_K^{-p})^\wedge$ is a positive distribution, i.e. $\langle \|\cdot\|_K^{-p},\varphi\rangle \ge 0$ for every non-negative test function $\varphi\in \mathcal S.$ By \cite[Theorem 6.16]{bib_kold_book}, for integers $k \in [1,n-1]$, $K$ embeds in $L_{-k}$ if and only if it is a $k$-intersection body (note that 1-intersection body is the usual intersection body of Lutwak). Also, every origin-symmetric convex body in $\R^n$ embeds in $L_{-n+3};$ see \cite[Corollary 4.9]{bib_kold_book}. Note that
the unit balls of finite dimensional subspaces of $L_p,\ 0<p\le 2$ embed in $L_{-p}$ for every $p\in (0,n).$ In particular, all ellipsoids embed in $L_{-p}.$
For details, we refer the reader to \cite[Section 6.3]{bib_kold_book}. We denote the class of bodies in $\R^n$ that embed in $L_{-p}$ by $L_{-p}^n.$

For a test function $\varphi:[0,\infty) \to \R$ the fractional derivative of order $q \in \C$ at $0$, is defined as the analytic continuation to $q \in \mathbb C$ of 
\[\varphi^{(q)}(0) = \frac 1{\Gamma(-q)} \int_0^\infty t^{-1-q} \varphi(t) dt.\]
If $\varphi$ is $m$-times continuously differentiable and integrable, then $\varphi^{(q)}(0)$ can be extended to an analytic function at $-1 < \Re(q) < m$.
Moreover, one has
\begin{multline}
	\varphi^{(q)}(0) = \frac 1{\Gamma(-q)} \int_0^1 t^{-1-q} \left( h(t)-h(0)- \cdots -h^{(m-1)}(0) \frac {t^{m-1}}{(m-1)!} \right) dt\\
	+ \frac 1{\Gamma(-q)} \int_1^\infty t^{-1-q} h(t) dt + \frac 1{\Gamma(-q)} \sum_{k=0}^{m-1} \frac{h^{(k)(0)}}{k!(k-q)}
\end{multline}
for $-1 < \Re(q) < m$, and the definition does not depend on the choice of $m$ as long as $\varphi$ is $m$-times continuously differentiable.

The Fourier transform of negative powers of the Euclidean norm $|\cdot|_2$ in $\R^n$ are given by
\begin{equation}
	\label{eq_transform_euclidean_norm}
(|\cdot|_2^{-q})^\wedge(\xi) = \frac{\pi^{n/2} \Gamma{\left(\frac{n-q}2\right)} 2^{n-q}}{\Gamma{(q/2)}} |\xi|_2^{q-n}
\end{equation}

for $q \in (0,n);$ see \cite{GS}.

\begin{definition} 
	For a function $f \in \S$ and $q \in (-n, \infty)$ the fractional Laplacian of the order $q/2$ is defined by
	\[\D f = (|\xi|^q \hat f(\xi))^\vee,\]
	where $\vee$ stands for the inverse Fourier transform.
\end{definition}
	
By \eqref{eq_transform_euclidean_norm}, for $q \in (0,n)$ we have the representation
\begin{equation}
	\label{eq_laplacian_convolution}
	\D[-q] f = c_{n,q} |\cdot|^{-n+q} * f,
\end{equation}
where
\[c_{n,q} = \frac{\pi ^{n/2} 2^{n-q} \Gamma \left(\frac{n-q}{2}\right)}{(2\pi)^n \Gamma \left(\frac q2\right)}.\]

\begin{lemma}
	\label{res_radon_laplacian_relation}
	For $q \in (-1, \infty)$ not an odd integer, the fractional Laplacian and the Radon transform are related by
	\[\Radq f = \Rad (\D f)(\xi, 0)\]
\end{lemma}
\begin{proof}
	By Lemma 2.23 in \cite{bib_kold_book}, if $q \in \mathbb C \setminus \{-1, -3, -5, \ldots \}$ and $t,\xi \in \R$,
	\[(|t|^{-1-q})^\wedge(\xi) = 2 \Gamma(-q) \cos(\pi q/2) |\xi|^q\]
	in the sense of analytic continuation of distributions.

	We use the well-known relation between the Radon and Fourier transforms; see \cite[Lemma 2.11]{bib_kold_book}.
	For every fixed $\xi\in S^{n-1}$ and any Schwartz test function $f$, the Fourier transform of the function $t \to \mathcal Rf(\xi,t)$ is the function $z\to \hat{f}(z\xi)$. Now, for $q \in (-1,0)$, we have
	\begin{align}
		\Radq f 
		&=\frac 1{\cos(\pi q/2)} \int_{-\infty}^\infty \frac 1{2 \Gamma(-q)} |t|^{-1-q} \Rad f(\xi, t) dt \\
		&= \frac 1{2\pi} \int_{-\infty}^\infty |z|^q \Rad f(\xi, \cdot)_t^\wedge(z) dz \\
		&= \frac 1{2\pi} \int_{-\infty}^\infty |z \xi|_2^q \hat f(z \xi) dz \\
		&= \frac 1{2\pi} \int_{-\infty}^\infty (\D f)^\wedge(z \xi) dz \\
		&= \frac 1{2\pi} \int_{-\infty}^\infty \Rad{(\D f)}(\xi, \cdot)_z^\wedge(z) dz \\
		&= \Rad{(\D f)}(\xi, 0).
	\end{align}
	By analytic continuation with respect to $q$ we obtain the equality for every $q \in (-1, \infty),$ not an odd integer.
\end{proof}

For the proof of the main results we need two lemmas.
\begin{lemma}
	\label{res_convolution}
	For positive measurable functions $f,g$, and measurable sets $A,B$,
	\[ \int_{A+B} f*g(x) dx \geq \int_A f(x) dx \int_B g(x) dx\]
\end{lemma}
\begin{proof}
	Since the support of $(f\chi_A)*(g\chi_B)$ is contained in $A+B$, we have
	\begin{align}
		\int_{A+B} f*g(x) dx 
		&\geq \int_{A+B} (f\chi_A)*(g\chi_B)(x) dx \\
		&\geq \int_{\R^n} (f\chi_A)*(g\chi_B)(x) dx \\
		&= \int_A f(x) dx \int_B g(x) dx \\
	\end{align}
\end{proof}

\begin{lemma}
	\label{res_Gamma_inequality}
	For $0 \leq \mu \leq \lambda$ we have
	\[\lambda^\mu \Gamma(\lambda - \mu) \geq \Gamma(\lambda).\]
\end{lemma}
\begin{proof}
	The equality holds for $\mu = 0$. Since the left-hand side is log-convex with respect to $\mu$, it suffices to show that
	its logarithmic derivative at $\mu = 0$ is non-negative.
	By Jensen's inequality,
	\begin{align}
	\left. \frac \partial{\partial \mu} \left[ \mu \log(\lambda) + \log \Gamma(\lambda - \mu) \right] \right|_{\mu=0}
		&= \log(\lambda) - \frac{\Gamma'(\lambda)}{\Gamma(\lambda)}\\
		&= \frac 1{\Gamma(\lambda)} \left( \int_0^\infty -\log(t/\lambda) t^{\lambda-1} e^{-t} dt \right) \\
		&\geq -\log \left( \frac 1{\lambda \Gamma(\lambda)} \int_0^\infty t^{\lambda} e^{-t} dt \right) \\
		&=0
	\end{align}
	where we used that $\Gamma(\lambda+1) = \lambda \Gamma(\lambda)$.
\end{proof}

\section{Proof of the Main Results}

\begin{proof}[Proof of Theorem \ref{res_fractional_example}]
	Consider $K$ and $f$ from Theorem \ref{res_example_KL}, and take $h = \D[-q] f$.
	By Lemma \ref{res_convolution} and equation \eqref{eq_laplacian_convolution},
	\begin{align}
		\int_{2K} h
		\label{eq_laplacian_bound}
		&\geq c_{n,q} \int_{\sqrt n\B} |x|_2^{-n+q} dx \int_K f dx \\
		&\geq c_{n,q} n\omega_n \frac{n^{q/2}}q c = c d_{n,q}.
	\end{align}
	where $d_{n,q} = \frac{2^{-q} n^{q/2} \Gamma \left(\frac{n-q}{2}\right)}{\Gamma \left(\frac n2\right) \Gamma \left(\frac q2+1\right)}$.

	By Lemma \ref{res_radon_laplacian_relation} and Theorem \ref{res_example_KL},
	\begin{align}
		\Radq h 
		\label{eq_radon_bound}
		&= \Radq{(\D[-q] f)} \\
		&= \Rad f(\xi, 0) \leq \frac c {\sqrt n}.
	\end{align}

	We consider $D= \vol{K}^{-1/n}K$ and
	\[g(x) = \left(\int_D h(\vol{2K}^{1/n} x) dx \right)^{-1} h(\vol{2K}^{1/n} x),\]
	so that $\int_D g = 1$ and $D$ has volume $1$.
	Combining \eqref{eq_laplacian_bound}, \eqref{eq_radon_bound} and applying Lemma \ref{res_Gamma_inequality}, the inequality $\vol{K} \leq c^n$ and Stirling's formula 
	\begin{align}
		\Radq g
		&= c \vol{2K}^{\frac {1+q}n} d_{n,q}^{-1} \Radq h \\
		&\leq c^q \frac 1{\sqrt n} \frac{2^q \Gamma \left(\frac n2\right) \Gamma \left(\frac q2+1\right)}{n^{q/2} \Gamma \left(\frac{n-q}{2}\right)} \\
		&\leq c^q \frac 1{\sqrt n} \Gamma \left(\frac q2+1\right) \\
		&\leq c^q \frac 1{\sqrt n} (q+1)^{\frac{q+1}2}.
	\end{align}
\end{proof}

%
%

\end{document}